\documentstyle[12pt, amsfonts]{amsart}
\begin{document}
\def\R{{\mathbb R}}
\def\Z{{\mathbb Z}}
\def\C{{\mathbb C}}
\newcommand{\trace}{\rm trace}
\newcommand{\Ex}{{\mathbb{E}}}
\newcommand{\Prob}{{\mathbb{P}}}
\newcommand{\E}{{\cal E}}
\newcommand{\F}{{\cal F}}
\newtheorem{df}{Definition}
\newtheorem{theorem}{Theorem}
\newtheorem{lemma}{Lemma}
\newtheorem{pr}{Proposition}
\newtheorem{co}{Corollary}
\def\n{\nu}
\def\sign{\mbox{ sign }}
\def\a{\alpha}
\def\N{{\mathbb N}}
\def\A{{\cal A}}
\def\L{{\cal L}}
\def\X{{\cal X}}
\def\F{{\cal F}}
\def\c{\bar{c}}
\def\v{\nu}
\def\d{\delta}
\def\diam{\mbox{\rm dim}}
\def\vol{\mbox{\rm Vol}}
\def\b{\beta}
\def\t{\theta}
\def\l{\lambda}
\def\e{\varepsilon}
\def\colon{{:}\;}
\def\pf{\noindent {\bf Proof :  \  }}
\def\endpf{ \begin{flushright}
$ \Box $ \\
\end{flushright}}

\title[Hyperplane inequality for measures]{A $\sqrt{n}$ estimate for measures of hyperplane sections of convex bodies}

\author{Alexander Koldobsky}

\address{Department of Mathematics\\ 
University of Missouri\\
Columbia, MO 65211}

\email{koldobskiya@@missouri.edu}

\begin{abstract}  The hyperplane (or slicing) problem asks whether there exists 
an absolute constant $C$ so that for any origin-symmetric convex body $K$ in $\R^n$
$$
|K|^{\frac {n-1}n} \le C \max_{\xi \in S^{n-1}} |K\cap \xi^\bot|,
$$
where  $\xi^\bot$ is the central hyperplane in $\R^n$ perpendicular to $\xi,$ and
$|K|$ stands for volume of proper dimension. The problem is still open, with the best-to-date estimate $C\sim n^{1/4}$ established
by Klartag, who slightly improved the previous estimate of Bourgain. It is much easier to get a weaker estimate with $C=\sqrt{n}.$

In this note we show that the $\sqrt{n}$ estimate holds for arbitrary measure in place of volume. Namely,
if $L$ is an origin-symmetric convex body in $\R^n$ and $\mu$ is a measure 
with non-negative even continuous density on $L,$ then
$$\mu(L)\ \le\ \sqrt{n} \frac n{n-1} c_n\max_{\xi \in S^{n-1}} 
\mu(L\cap \xi^\bot)\ |L|^{1/n} \ ,$$
where  $c_n= \left|B_2^n\right|^{\frac{n-1}n}/ \left|B_2^{n-1}\right| < 1,$
and $B_2^n$ is the unit Euclidean ball in $\R^n.$ We deduce this inequality from a stability
result for intersection bodies.

\end{abstract}  
\maketitle

\section{Introduction}
The hyperplane (or slicing) problem  \cite{Bo1, Bo2, Ba, MP} asks whether there exists 
an absolute constant $C$ so that for any origin-symmetric convex body $K$ in $\R^n$
\begin{equation} \label{hyper}
|K|^{\frac {n-1}n} \le C \max_{\xi \in S^{n-1}} |K\cap \xi^\bot|,
\end{equation}
where  $\xi^\bot$ is the central hyperplane in $\R^n$ perpendicular to $\xi,$ and
$|K|$ stands for volume of proper dimension.
The problem is still open, with the best-to-date estimate $C\sim n^{1/4}$ established
by Klartag \cite{Kl}, who slightly improved the previous estimate of Bourgain \cite{Bo3}.
We refer the reader to [BGVV] for the history and 
current state of the problem.

In the case where $K$ is an intersection body (see definition and properties below), 
the inequality (\ref{hyper}) can be proved with 
the best possible constant (\cite[p. 374]{G2}):
\begin{equation}\label{hyper-inter}
|K|^{\frac {n-1}n} \le \frac{\left|B_2^n\right|^{\frac{n-1}n}}{\left|B_2^{n-1}\right|} 
\max_{\xi \in S^{n-1}} |K\cap \xi^\bot|,
\end{equation}
with equality when $K=B_2^n$ is the unit Euclidean ball.  Here $|B_2^n|= \pi^{n/2}/\Gamma(1+n/2)$
is the volume of $B_2^n.$ Throughout the paper, we denote the constant in (\ref{hyper-inter}) by
$$c_n= \frac{\left|B_2^n\right|^{\frac{n-1}n}}{\left|B_2^{n-1}\right|} .$$
Note that $c_n<1$ for every $n\in \N;$ this is an easy consequence of the log-convexity 
of the $\Gamma$-function. 

It was proved in \cite{K3}  that inequality (\ref{hyper}) holds for intersection bodies
with arbitrary measure in place of volume. Let $f$ be an even continuous non-negative
function on $\R^n,$ and denote by $\mu$ the measure on $\R^n$ with density $f$. 
For every closed bounded set $B\subset \R^n$ define
$$\mu(B)=\int\limits_B f(x)\ dx.$$  Suppose that $K$ is an intersection body in $\R^n.$ Then,
as proved in \cite[Theorem 1]{K3} (see also a remark at the end of the paper \cite{K3}),
\begin{equation} \label{arbmeas}
\mu(K) \le \frac n{n-1} c_n \max_{\xi \in S^{n-1}} \mu(K\cap \xi^\bot)\ |K|^{1/n}.
\end{equation}
The constant in the latter inequality is the best possible. 

This note was motivated by a question of whether one can  remove the assumption that
$K$ is an intersection body and  prove the inequality (\ref{arbmeas}) for all origin-symmetric convex
bodies, perhaps at the expense of a greater constant in the right-hand side. One would like this extra 
constant to be independent of the body or measure. In this note we prove the following inequality.

\begin{theorem}\label{main} Let $L$ be an origin-symmetric convex body in $\R^n,$ and 
let $\mu$ be a measure with even continuous non-negative density on $L.$ Then
\begin{equation} \label{sqrtn}
\mu(L)\ \le\ \sqrt{n} \frac n{n-1} c_n\max_{\xi \in S^{n-1}} 
\mu(L\cap \xi^\bot)\ |L|^{1/n}.
\end{equation}
\end{theorem}

In the case of volume, the estimate (\ref{hyper}) with $C=\sqrt{n}$ can be proved relatively
easy (see \cite[p. 96]{MP} or \cite[Theorem 8.2.13]{G2}), and it is not optimal, as mentioned above.
The author does not know whether the estimate (\ref{sqrtn}) is optimal for arbitrary measures.

\section{Proof of Theorem \ref{main}}

We need several definitions and facts.
A closed bounded set $K$ in $\R^n$ is called a {\it star body}  if 
every straight line passing through the origin crosses the boundary of $K$ 
at exactly two points different from the origin, the origin is an interior point of $K,$
and the {\it Minkowski functional} 
of $K$ defined by 
$$\|x\|_K = \min\{a\ge 0:\ x\in aK\}$$
is a continuous function on $\R^n.$ 

The {\it radial function} of a star body $K$ is defined by
$$\rho_K(x) = \|x\|_K^{-1}, \qquad x\in \R^n.$$
If $x\in S^{n-1}$ then $\rho_K(x)$ is the radius of $K$ in the
direction of $x.$

If $\mu$ is a measure on $K$ with even continuous density $f$, then 
\begin{equation} \label{polar-measure}
\mu(K) = \int_K f(x)\ dx = \int\limits_{S^{n-1}}\left(\int\limits_0^{\|\theta\|^{-1}_K} r^{n-1} f(r\theta)\ dr\right) d\theta.
\end{equation}
Putting $f=1$, one gets
\begin{equation} \label{polar-volume}
|K|
=\frac{1}{n} \int_{S^{n-1}} \rho_K^n(\theta) d\theta=
\frac{1}{n} \int_{S^{n-1}} \|\theta\|_K^{-n} d\theta.
\end{equation}

The {\it spherical Radon transform} 
$R:C(S^{n-1})\mapsto C(S^{n-1})$  
is a linear operator defined by
$$Rf(\xi)=\int_{S^{n-1}\cap \xi^\bot} f(x)\ dx,\quad \xi\in S^{n-1}$$
for every function $f\in C(S^{n-1}).$

The polar formulas (\ref{polar-measure}) and  (\ref{polar-volume}), applied to a hyperplane section of $K$, express 
volume of such a section in terms of the spherical Radon transform:

$$\mu(K\cap \xi^\bot) = \int_{K\cap \xi^\bot} f =  
\int_{S^{n-1}\cap \xi^\bot} \left(\int_0^{\|\theta\|_K^{-1}} r^{n-2}f(r\theta)\ dr \right)d\theta$$
\begin{equation} \label{measure=spherradon}
=R\left(\int_0^{\|\cdot\|_K^{-1}} r^{n-2}f(r\ \cdot)\ dr \right)(\xi).
\end{equation}
and
\begin{equation} \label{volume=spherradon}
|K\cap \xi^\bot| = \frac{1}{n-1} \int_{S^{n-1}\cap \xi^\bot} \|\theta\|_K^{-n+1}d\theta =
\frac{1}{n-1} R(\|\cdot\|_K^{-n+1})(\xi).
\end{equation}

The spherical Radon 
transform is self-dual (see \cite[Lemma 1.3.3]{Gr}), namely,
for any functions $f,g\in C(S^{n-1})$
\begin{equation} \label{selfdual}
\int_{S^{n-1}} Rf(\xi)\ g(\xi)\ d\xi = \int_{S^{n-1}} f(\xi)\ Rg(\xi)\ d\xi.
\end{equation}
Using self-duality, one can extend the spherical Radon transform to measures. 
Let $\mu$ be a finite Borel measure on $S^{n-1}.$
We define the spherical Radon transform of $\mu$ as a functional $R\mu$ on
the space $C(S^{n-1})$ acting by
$$(R\mu,f)= (\mu, Rf)= \int_{S^{n-1}} Rf(x) d\mu(x).$$
By Riesz's characterization of continuous linear functionals on the
space $C(S^{n-1})$, 
$R\mu$ is also a finite Borel measure on $S^{n-1}.$ If $\mu$ has 
continuous density $g,$ then by (\ref{selfdual}) the 
Radon transform of $\mu$ has density $Rg.$

The class of intersection bodies was introduced by Lutwak \cite{L}.
Let $K, L$ be origin-symmetric star bodies in $\R^n.$ We say that $K$ is the 
intersection body of $L$ if the radius of $K$ in every direction is 
equal to the $(n-1)$-dimensional volume of the section of $L$ by the central
hyperplane orthogonal to this direction, i.e. for every $\xi\in S^{n-1},$
\begin{equation} \label{intbodyofstar}
\rho_K(\xi)= \|\xi\|_K^{-1} = |L\cap \xi^\bot|.
\end{equation} 
All bodies $K$ that appear as intersection bodies of different star bodies
form {\it the class of intersection bodies of star bodies}. 

Note that the right-hand
side of (\ref{intbodyofstar}) can be written in terms of the spherical Radon transform using (\ref{volume=spherradon}):
$$\|\xi\|_K^{-1}= \frac{1}{n-1} \int_{S^{n-1}\cap \xi^\bot} \|\theta\|_L^{-n+1} d\theta=
\frac{1}{n-1} R(\|\cdot\|_L^{-n+1})(\xi).$$
It means that a star body $K$ is 
the intersection body of a star body if and only if the function $\|\cdot\|_K^{-1}$
is the spherical Radon transform of a continuous positive function on $S^{n-1}.$
This allows to introduce a more general class of bodies. A star body
$K$ in $\R^n$ is called an {\it intersection body} 
if there exists a finite Borel measure \index{intersection body}
$\mu$ on the sphere $S^{n-1}$ so that $\|\cdot\|_K^{-1}= R\mu$ as functionals on 
$C(S^{n-1}),$ i.e. for every continuous function $f$ on $S^{n-1},$
\begin{equation} \label{defintbody}
\int_{S^{n-1}} \|x\|_K^{-1} f(x)\ dx = \int_{S^{n-1}} Rf(x)\ d\mu(x).
\end{equation}

We refer the reader to the books \cite{G2, K2}
for more information about intersection bodies and their applications. Let  us just say that
intersection bodies played a crucial role in the solution of the Busemann-Petty problem.
The class of intersection bodies is rather rich. For example, every origin-symmetric convex 
body in $\R^3$ and $\R^4$ is an intersection body \cite{G1,  Z}. The unit ball of any finite
dimensional subspace of $L_p,\ 0<p\le 2$ is an intersection body, in particular every polar
projection body is an intersection body \cite{K1}. 
 
We deduce Theorem 1 from the following stability result for intersection bodies.  
\begin{theorem}\label{stab}
Let $K$ be an intersection body in $\R^n,$ let $f$
be an even continuous function on $K,$ $f\ge 1$ everywhere on $K,$ and let $\e>0.$ Suppose that 
\begin{equation}\label{comp1}
\int_{K\cap \xi^\bot} f \ \le\ |K\cap \xi^\bot| +\e,\qquad \forall \xi\in S^{n-1},
\end{equation}
then
\begin{equation}\label{comp2}
\int_K f\ \le\ |K| + \frac {n}{n-1}\ c_n\ |K|^{1/n}\e.
\end{equation}
\end{theorem}

\pf First, we use the polar formulas (\ref{measure=spherradon}) and (\ref{volume=spherradon}) to write
the condition (\ref{comp1}) in terms of the spherical Radon transform:
$$R\left(\int_0^{\|\cdot\|_K^{-1}} r^{n-2}f(r\ \cdot)\ dr \right)(\xi) \le \frac{1}{n-1} R(\|\cdot\|_K^{-n+1})(\xi) + \e.$$
Let $\mu$ be the measure on $S^{n-1}$ corresponding to $K$ by the definition of an intersection body (\ref{defintbody}).
Integrating both sides of the latter inequality over $S^{n-1}$ with the measure $\mu$ and using (\ref{defintbody}),
we get
$$\int_{S^{n-1}} \|\theta\|_K^{-1} \left(\int_0^{\|\theta\|_K^{-1}} r^{n-2}f(r\theta)\ dr \right)d\theta $$
\begin{equation} \label{eq11}
\le  \frac{1}{n-1} \int_{S^{n-1}} \|\theta\|_K^{-n}\ d\theta + \e \int_{S^{n-1}} d\mu(\xi).
\end{equation}
Recall (\ref{polar-measure}), (\ref{polar-volume}) and the assumption that $f\ge 1.$ We write the integral in the left-hand side 
of (\ref{eq11}) as
$$\int_{S^{n-1}} \|\theta\|_K^{-1} \left(\int_0^{\|\theta\|_K^{-1}} r^{n-2}f(r\theta)\ dr \right)d\theta $$
$$= \int_{S^{n-1}}  \left(\int_0^{\|\theta\|_K^{-1}} r^{n-1}f(r\theta)\ dr \right)d\theta$$ 
$$+ \int_{S^{n-1}} \left(\int_0^{\|\theta\|_K^{-1}} (\|\theta\|_K^{-1} - r)  r^{n-2}f(r\theta)\ dr \right)d\theta$$
$$\ge \int_K f + \int_{S^{n-1}} \left(\int_0^{\|\theta\|_K^{-1}} (\|\theta\|_K^{-1} - r)  r^{n-2}\ dr \right)d\theta$$
\begin{equation} \label{eq22}
=\int_K f + \frac 1{(n-1)n} \int_{S^{n-1}} \|\theta\|_K^{-n}\ d\theta = \int_K f + \frac1{n-1} |K|.
\end{equation}

Let us estimate the second term in the right-hand side of (\ref{eq11}) by adding the Radon transform of  the unit constant
 function under the integral ($R1(\xi)=\left|S^{n-2}\right|$ for every $\xi \in S^{n-1}$), 
 using again the fact that $\|\cdot\|_K^{-1}=R\mu$ and then applying H\"older's  inequality:
$$\e \int_{S^{n-1}} d\mu(\xi) = \frac{\e}{\left|S^{n-2}\right|} \int_{S^{n-1}} R1(\xi)\ d\mu(\xi)$$
$$=\frac{\e}{\left| S^{n-2} \right| } \int_{S^{n-1}} \|\theta\|_K^{-1}\ d\theta $$
$$ \le  \frac{\e}{\left|S^{n-2}\right|} \left|S^{n-1}\right|^{\frac{n-1}n} \left(\int_{S^{n-1}} \|\theta\|_K^{-n}\ d\theta\right)^{\frac1n}$$
\begin{equation}\label{eq33}
=  \frac{\e}{\left|S^{n-2}\right|} \left|S^{n-1}\right|^{\frac{n-1}n} n^{1/n}|K|^{1/n}= \frac n{n-1} c_n |K|^{1/n} \e.
\end{equation}
In the last step we used $|S^{n-1}|=n|B_2^n|, |S^{n-2}|=(n-1)|B_2^{n-1}|.$ Combining (\ref{eq11}), (\ref{eq22}),(\ref{eq33}) we get
$$\int_K f + \frac 1{n-1} |K| \le \frac n{n-1} |K| + \frac n{n-1} c_n |K|^{1/n} \e. \qed$$
\bigbreak
Now we prove our main result.
\smallbreak
\noindent {\bf Proof of Theorem \ref{main}: }
Let $g$ be the density of the measure $\mu,$ so $g$ is an even non-negative continuous function
on $L.$ By John's theorem \cite{J}, there exists an origin-symmetric ellipsoid $K$ such that
$$\frac 1{\sqrt{n}} K \subset L \subset K.$$ 
The ellipsoid $K$ is an intersection body (see for example \cite[Corollary 8.1.7]{G2}). 
Let $f= \chi_K + g \chi_L,$ where $\chi_K,\ \chi_L$ are the indicator functions of $K$ and $L.$ 
Clearly, $f\ge 1$ everywhere on $K.$ Put 
$$\e=\max_{\xi\in S^{n-1}} \left(\int_{K\cap \xi^\bot} f - |K\cap \xi^\bot| \right)= \max_{\xi\in S^{n-1}} \int_{L\cap \xi^\bot} g$$
and apply Theorem \ref{stab} to $f,K,\e$ (the function $f$ is not necessarily continuous on $K,$ 
but the result holds by a simple approximation argument). We get
$$\mu(L)= \int_L g = \int_K f  -\ |K|$$
$$ \le \frac n{n-1} c_n |K|^{1/n}\max_{\xi\in S^{n-1}} \int_{L\cap \xi^\bot} g$$
$$ \le \sqrt{n}\ \frac n{n-1} c_n |L|^{1/n}\max_{\xi\in S^{n-1}} \mu(L\cap \xi^\bot),$$
because $|K|^{1/n}\le \sqrt{n}\ |L|^{1/n}.$ \qed

\bigbreak
{\bf Acknowledgement.} I wish to thank the US National Science Foundation for support through 
grant DMS-1265155.

\end{document}